\documentclass[reqno]{article}
\usepackage{amssymb, amsmath, amsthm, amsfonts, amscd}
\usepackage[english]{babel}

\usepackage{graphicx}

\usepackage{epstopdf}

\usepackage{breqn}
\usepackage{mathtools}
\usepackage{color}
\usepackage{hyperref}
\hypersetup{
    colorlinks=true,
    citecolor=red,
    filecolor=black,
    linkcolor=blue,
    urlcolor=black
}

\chardef\bslash=`\\ 

\newtheorem{thm}{Theorem}[section]
\newtheorem{cor}[thm]{Corollary}
\newtheorem{lem}[thm]{Lemma}
\newtheorem{prop}[thm]{Proposition}

\newtheorem{defn}{Definition}[section]
\newtheorem{notation}{Notation}[section]

\newcommand{\N}{\mathbb{N}}
\newcommand{\NN}{\mathcal{N}}
\newcommand{\EE}{\mathcal{E}}

\newcommand{\PP}{\mathcal{P}}

\newcommand{\VV}{\mathcal{V}}

\newcommand{\Z}{\mathbb{Z}}
\newcommand{\Q}{\mathbb{Q}}
\newcommand{\R}{\mathbb{R}}
\newcommand{\C}{\mathbb{C}}
\newcommand{\Sp}{\mathbb{S}}
\newcommand{\T}{\mathbb{T}}
\newcommand{\TT}{\mathcal{T}}

\def\a{\alpha }

\def\l{\lambda }

\def\p{\pi}

\def\s{\sigma}

\def\z{\zeta}
\def\w{\omega}
\def\e{\varepsilon}
\def\f{\varphi}
\def\.{\cdot }
\def\ra{\rightarrow}
\def\hra{\hookrightarrow}

\title{
\textsc{\textbf{Continuous spectrum or measurable reducibility for
quasiperiodic cocycles in $\T ^{d} \times SU(2)$}}\\
\author{Nikolaos Karaliolios \footnote{IME, UFF, Niteroi, RJ, Brazil.
email: nikos@math.uff.br}
}}

\begin{document}

\maketitle

\begin{abstract}
We continue our study of the local theory for quasiperiodic
cocycles in $\T ^{d} \times G$, where $G=SU(2)$, over a rotation
satisfying a Diophantine condition and satisfying a closeness-to-constants
condition, by proving a dichotomy between measurable reducibility (and therefore pure point spectrum), and
purely continuous spectrum in the space orthogonal to
$L^{2}(\T ^{d}) \hra L^{2}(\T ^{d} \times G)$. Subsequently,
we describe the equivalence classes of cocycles under smooth
conjugacy, as a function of the parameters defining their
K.A.M. normal form.
Finally, we derive a complete classification of the dynamics of one-frequency ($d=1$) cocycles
over a Recurrent Diophantine rotation.

All theorems will be stated sharply in terms of the number of
frequencies $d$, but in the proofs we will always assume $d=1$, for
simplicity in expression and notation.
\end{abstract}  

\tableofcontents

\nocite{ParryTopicsErgTh}
\nocite{KatokThouv2006}
\nocite{FayadKatok2004}

\section{Introduction and statement of the results}
%
%
This article continues the work taken up by the author in \cite{NKPhD}, \cite{NKRigidity} and
\cite{NKInvDist}. This work used and developed the techniques of \cite{KrikAst}, \cite{Krik2001},
\cite{El2002a} and \cite{Fra2004} among others, and this
part was motivated by \cite{deAldecoa13}, where conditions for the existence
of absolutely continuous spectrum are given.

We think that the K.A.M.-theoretical, or local, part of the theory
of such dynamical systems from the topological point of view is
concluded with the following theorem and with the classification
and path connectedness theorems that we state later on. On the
other hand, the metric abundance of reducible cocycles in the
$C^{\infty}$ category is an open question (see \cite{KrikAst} for
the proof of this theorem in the analytic category).
\begin{thm} \label{thm spectral dichotomy}
Let $\a \in \T ^{d}$ satisfy a Diophantine condition $DC(\gamma , \tau )$. Then, there exist
$\e >0 $ and $s_{0} \in \N ^{*}$, depending on $d , \gamma , \tau $, such that every cocycle
$(\a , A.e^{F(\. )}) $ with $A \in G = SU(2)$,
$\| F \|_{0} < \e $ and $\| F \|_{s_{0}} < 1 $
satisfies the following dichotomy:
\begin{enumerate}
\item \label{meas. red} either the cocycle is measurably reducible and
therefore has pure-point spectrum,
\item \label{cont spec} or the spectrum of the associated Koopman operator is purely continuous
in the space orthogonal to $L^{2}(\T ^{d}) \hra L^{2}(\T ^{d} \times G,
\C)$. Such cocycles are weak mixing in the fibers, and they are not strong
mixing.
\end{enumerate}
\end{thm}
The neighborhood of constants described in the statement of the
theorem will be
referred to as the \textit{K.A.M.} or the \textit{local regime} and
denoted by $\NN$. Conjugations of the size of those produced by the
K.A.M. scheme when the product thus constructed converges, or of those
who reduce a given cocycle to its K.A.M. normal form, satisfy
a condition of the type $|Y|_{0}< C $
and $|Y|_{s_{0} - \xi \gamma }< 1$, for some constants
$\epsilon \leq C <1 $ and $ 1 \leq \xi < s_{0} / \gamma $
depending on $\gamma , \tau$ and $d$. Conjugations satisfying
such estimates will be referred to as \textit{close-to-the-identity conjugations}. They form a contractible
set in $C^{\infty}(\T ^{d},G)$ denoted by $\VV$.

We make the following remarks. Firstly, measurable reducibility, as
in item \ref{meas. red}, is equivalent to a condition which is
explicit in terms of the output of the K.A.M. scheme with the parameters of theorem
\ref{thm spectral dichotomy}, (cf. \cite{NKInvDist} or the discussion below).
The question of differentiable
rigidity of measurable reducibility was investigated in \cite{NKRigidity},
where we proved that measurable reducibility to a full measure
set of constants $DC_{\a} \subset G$
implies in fact smooth reducibility (see thm \ref{thm dif rig}
below). On the other hand,
in \cite{NKInvDist} it was shown that for constants in
the generic set $\mathcal{L}_{a} = G \setminus DC_{\a}$,
reducibility is not rigid. Moreover, the construction of the
measurable transfer function under the relevant condition shows that its form is quite special.
The assumption that a measurable conjugation reduces a smooth cocycle imposes some
constraints on the former, which consequently will not be very wild.

Secondly, the non-existence of a measurable conjugation reducing the
coycle to a constant, as in item \ref{cont spec}, is equivalent to
the complementary condition to the one in item \ref{meas. red}, and
thus equally explicit. In the complementary space of the one bearing
the purely continuous spectrum (i.e. in $L^{2}(\T ^{d})$, the Kronecker factor of the cocycle),
the spectrum is pure-point, because of the quasiperiodic
dynamics within the torus.

The proof is based on the use of the K.A.M. normal form, introduced in \cite{NKInvDist},
in order to prove that the existence of
an eigenfunction of the Koopman operator associated to a given cocycle implies the
condition for measurable reducibility, provided that the eigenfunction depends
non-trivially on the variable in the fibers. The K.A.M. normal form, though not essential
to the proof, greatly simplifies the calculations and elucidates the geometry of the
problem. The existence of the K.A.M. normal form is a corollary of
the almost reducibility theorem, first obtained in \cite{El2002a},
and its function in our study, put informally, is to separate
the close-to-the-identity part of the conjugation (the one in $\VV$) from the
far-from-the-identity part, both constructed almost exactly as in \cite{El2002a},
and keep
the second which contains all the interesting information. The purpose of
the subsequent analysis is to establish the (necessary and sufficient, as
it turns out) conditions under which some rearrangement of these
far-from-the-identity conjugations converges in some function space. The
(necessary and sufficient) condition for measurable reducibility
that we referred to above, in this context, is that the angles
$\theta _{i}$ between the successive far-from-the-identity
conjugations at the steps $i$ and $i+1$ of the scheme be square
summable. Still informally, $\theta _{i}$ estimates the
commutator appearing in the following sequence of operations: solve an equation
in the coordinates of the $i$-th step of the scheme, make one step
of scheme, solve the same equation with greater precision,
undo the change of coordinates and check the compatibility of
the expressions. We remark that cohomology in  $C^{\infty}(\T , G)$
over the rotation $x \mapsto x+\a $ seems to demand only one
step of the scheme (see the proof of thms 1.1 and 1.2 in
\cite{NKInvDist} where the rearrangement of the far-from-the-identity
conjugations is described) for the estimate to be effective, while
cohomology in $C^{\infty}$ or in $L^{2}(\T ^{d}\times G , \C )$
over any given cocycle seems to demand two steps of the scheme (see
e.g. the proof of thm \ref{cont spec} herein).

The proof that non-reducible cocycles are not strong mixing uses the
fact that a dynamical system is stong mixing iff the iterates of the
Koopman operator associated to it converge to $0$ in the weak operator
topology. This property is incompatible with rigidity, in the sense
that for every cocycle $(\a , A(\. ))$ in the local regime, there
exists a sequence of iterates $m_{j} \ra \infty$ such that
\begin{equation*}
(\a , A(\. ))^{m_{j}} \overset{C^{\infty}}\ra (0 , Id)
\end{equation*}
%

\bigskip

Combining the results obtained in the literature mentioned in the beginning of
this section with thm \ref{thm spectral dichotomy}, we obtain the following picture for the dynamics of cocycles in
$\T ^{d} \times G$, close to a constant, supposing that we have run the K.A.M. scheme
(\cite{KrikAst}, \cite{El2002a}, \cite{NKRigidity}, \cite{NKInvDist}, \cite{NKPhD}).

\begin{thm} [\cite{KrikAst},\cite{El2002a},\cite{NKInvDist}]
Every cocycle in the K.A.M. regime, $\NN $, is almost reducible, and conjugate
to a cocycle in K.A.M. normal form.
\end{thm}
The next theorems concern the reducibility of cocycles, via transfer
functions of regularity $L^{2}$, $C^{\infty}$, or intermediate
Sobolev regularity $H^{s}, 0<s<\infty$.
\begin{thm} [\cite{NKInvDist}] \label{thm meas reduc}
A given cocycle is measurably reducible iff the angles of successive
far-from-the-identity conjugations are summable in $\ell ^{2}$.
We can then construct a sequence of $C^{\infty}$ smooth
conjugations converging in $L^{2}$ toward a reducing conjugation.
\end{thm}
The following is a corollary of the proof of the previous theorem.
\begin{cor} [\cite{NKInvDist}]
If the angles are in fact summable in higher regularity $h ^{\s} , 0<\s\leq \infty$ (i.e. if $\{ N_{n_{i}}^{\s} \theta _{i} \} \in \ell
^{2} $), then
the cocycle is $H^{\s}$-reducible and the conjugation constructed as above converges
in $H^{\s}$.
\end{cor}
A particular case of the corollary is that of the occurrence of a finite number of
far-from-the-identity conjugations, where summability in every
$H^{\s }$ is trivial. This always occurs when the cocycle is reducible
to a constant which is "more Diophantine than the frequency"
(cf. \cite{NKRigidity}).
\begin{thm} [\cite{NKInvDist}]
For any $\s, 0\leq \s \leq \infty $, cocycles that are reducible in any given regularity
$H^{\s}$, and not any higher, are dense in $\NN $. Reducibility in any given regularity
is an $F_{\s }$ condition.
\end{thm}
The following theorem states that for "a generic" reducible
cocycle measurable conjugation is not rigid, while "for almost every"
one, it is.
\begin{thm} [\cite{NKRigidity}, \cite{NKInvDist}] \label{thm dif rig}
Reducibility in regularity $0\leq \s < \infty$ can only
occur when the constant cocycle $(\a , A)$ to which
the cocycle is reduced is Liouville with respect to $\a$, i.e.
when $A \in \mathcal{L} _{\a}$. If the constant cocycle is
Diophantine with respect to $\a$, i.e. if
$A \in DC_{\a } = G \setminus \mathcal{L}_{\a}$,
then measurable reducibility implies smooth reducibility.
\end{thm}

The following theorems, starting with thm \ref{thm spectral dichotomy} concern non-reducible cocycles.
Herein, we prove that all cocycles that are not measurably
reducible are weak mixing in the fibers, and therefore uniquely
ergodic in the whole space for the product of the Haar measures.
Non-reducibility is a $G_{\delta }$ condition, and is also dense in
$\NN $. This theorem strengthens H. Eliasson's theorem
in two ways. The condition we obtain is striclty looser than the one
given by H. Eliasson, as well as being optimal. Moreover, we prove
that the cocycles satisfying this condition are not only uniquely
ergodic, but in fact weak mixing. Generic cocycles have a stronger
property of ergodicity.

\begin{thm} [\cite{NKInvDist}]
Distributional unique ergodicity is also equivalent to a
$G_{\delta }$-dense condition, which is stricter than the one for weak mixing.
\end{thm}
We refer the reader unfamiliar with the notions to \cite{NKInvDist}
and to the bibliography therein for the definition of DUE and the 
related concepts.

Distributional unique ergodicity is equivalent to an
explicit condition on the asymptotic repartition of the angles
between successive far-from-the-identity conjugations. Relaxation
of this condition in a controlled way gives the following theorem.
\begin{thm} [\cite{NKInvDist}]
In the border between unigue ergodicity in the space of distributions and in the
classical sense, countably infinite invariant distributions of arbirtarily high orders are created.
\end{thm}
We have also proved that DUE cocycles are not Cohomologically
Stable, since the following theorem holds. The Diophantine condition
of the theorem is stricter than $DC_{\a}$, see paragraph
\ref{Arithmetics, continued fraction expansion}.
\begin{thm} [\cite{NKInvDist}] \label{thm cohom stab impl red}
A cocycle in $\NN $ is cohomologically stable iff it is
$C^{\infty}$ reducible to a constant cocycle $(\a, A)$ inducing a
Diophantine rotation on its invariant tori
$\TT _{(\a, A)} \approx \T ^{d} \times \Sp ^{1}
\hra \T ^{d} \times G $.
\end{thm}

The following theorems concern the topology of the different
conjugacy classes, and the way they lie in $\NN $.
\begin{thm} \label{thm connect red to any in N}
Given any cocycle $(\a ,A^{\prime} (\. )) \in \NN$, one can construct
a continuous path $[0,1] \ra \NN $ such that $(\a ,A_{0} (\. ))= (\a, Id)$,
for all $t \in [0,1)$ the cocycle $(\a ,A_{t} (\. ))$ is $C^{\infty}$
reducible, and $(\a ,A_{1} (\. ))$ is the K.A.M. normal form
of $(\a ,A^{\prime} (\. ))$.
\end{thm}
The path is in fact  piecewise $C^{\infty}$. If we allow the path
to exit the K.A.M. regime, we can obtain more.
\begin{thm} \label{thm connect red to any outside N}
Given any cocycle $(\a ,A^{\prime} (\. )) \in \NN$, one can construct
a continuous path $[0,1] \ra SW^{\infty}_{\a} (\T ^{d},G) $ such that $(\a ,A_{0} (\. ))= (\a, Id)$,
for all $t \in [0,1)$ the cocycle $(\a ,A_{t} (\. ))$ is $C^{\infty}$
conjugate to $(\a, Id)$, and $(\a ,A_{1} (\. ))$ is the
K.A.M. normal form of $(\a ,A^{\prime} (\. ))$.
\end{thm}
The path $\gamma \colon [0,1) \ra C^{\infty}(\T ^{d}, G) $
acting on $(\a, Id)$ and producing
\begin{equation*}
(\a ,A_{t} (\. )) = Conj_{\gamma (t)}(\a, Id)
\end{equation*}
is continuous (in fact piecewise $C^{\infty}$), and
degenerates in a prescribed way when the target cocycle is not
$C^{\infty }$ reducible.
At time $t= 1^{-}$ it may exit $C^{\infty}$ into a function space of lower regularity,
or a space of distributions\footnote{In fact the limit is always well
defined in $H^{-d/2}(\T ^{d}, G)$.}.
The path in the conjugacy space may exit $\VV$, and in general will do
so. The path in $SW^{\infty}_{\a}$ will consequently exit $\NN$ and reenter,
in general an infinite number of times.
Since the space of conjugations taking a cocycles to their respective
normal forms is contractible (it is the space $\VV$), we immediately obtain the following corollary.
\begin{cor}
The target of the path can be the cocycle $(\a ,A^{\prime} (\. ))$
itself, while the same as in theorem \ref{thm connect red to any in N}
(resp. thm \ref{thm connect red to any outside N}) holds for all times
$t \in [0,1)$.
\end{cor}
We can in fact obtain the following, stronger theorem, which
establishes a way in which the topology of any two classes share some
properties.
\begin{thm} \label{thm connect any to any in N}
Given any two cocycles
$(\a ,A_{0} (\. ))$ and $(\a ,A' (\. ))$, in $\NN $,
one can construct a continuous path $[0,1] \ra \NN $, and such that:
\begin{enumerate}
\item for every $t \in [0,1)$ the normal form of the cocycle
$(\a ,A_{t} (\. ))$ has the same as tail as that of $(\a ,A_{0} (\. ))$,
and thus the same dynamical properties.
\item $(\a ,A_{1} (\. ))$ is the K.A.M. normal form of
$(\a ,A' (\. ))$, and thus conjugate to it via
a conjugation in $\VV$.
\end{enumerate}
As before, the target cocycle can be the cocycle $(\a ,A' (\. ))$.
\end{thm}
Informally, we can connect any two types of dynamical behavior
without exiting $\NN$. The following theorem says that if we allow
the path to exit $\NN$, we can do the same thing with conjugacy classes.
\begin{thm} \label{thm connect any to any outside N}
Given any two cocycles
$(\a ,A_{0} (\. ))$ and $(\a ,A' (\. ))$, in $\NN $, one can construct
a continuous path $[0,1] \ra SW^{\infty}_{\a} $ such that:
\begin{enumerate}
\item for every $t \in [0,1)$, the cocycle $(\a ,A_{t} (\. ))$ is
conjugate to $(\a ,A_{0} (\. ))$.
\item $(\a ,A_{1} (\. ))$ is the K.A.M. normal form of
$(\a ,A' (\. ))$, and thus conjugate to it via
a conjugation in $\VV$.
\end{enumerate}
As before, the target cocycle can be the cocycle $(\a ,A' (\. ))$.
\end{thm}
Again, a path in the space of conjugations acting on
$(\a ,A_{0} (\. ))$ is constructed with the same properties as in
theorem \ref{thm connect red to any outside N}.

These two last theorems illustrate the necessity of a transversality
condition for obtaining a full-measure reducibility theorem for one-parameter
families of cocycles as in \cite{Krik99} or \cite{KrikAst}.
On the other hand, the K.A.M. normal form should be expected
to depend badly on parameters along a generic family, so our
approach is not expected to be well adapted to the metric point of view.

Finally, we give a satisfactory classification of conjugation
classes.
\begin{thm} \label{thm classif of classes}
Two cocycles in K.A.M. normal form are $C^{\infty}$
conjugate iff the parameters defining the normal forms satisfy
the following properties.
\begin{enumerate}
\item The resonant steps $n_{i}^{j}$, $j=1,2$, are the same, except for
a finite number.
\item The angles between successive conjugations
\begin{equation*}
\theta ^{j} _{i} = \arctan \frac{|\hat{F}^{j}_{i}(k^{j}_{i}) |}
{|\epsilon^{j} _{i}|} , j=1,2
\end{equation*}
are equal up to $O(N_{i-1} ^{-\infty})$.
\item The "rotation numbers"
\begin{equation*}
a^{j}_{i } = k^{j}_{i } \a  \mod \Z 
 \simeq \sqrt{|\hat{F}^{j}_{i-1}(k^{j}_{i-1}) |^{2} + ( \epsilon^{j} _{i-1}) ^{2}}, j=1,2
\end{equation*}
satisfy
\begin{equation*}
a^{1}_{i } - a^{2}_{i } = \tilde{k}_{i} \a \mod \Z
\end{equation*}
and
\begin{equation*}
k^{1}_{i } -k^{2}_{i } = \tilde{k}_{i}
\end{equation*}
Here, $\tilde{k}_{i} \in \Z  $ has to
be such that $N_{i-1} <  k^{j}_{i}  \leq N_{i}, j=1,2$,
and $\tilde{k}_{i} \in \Z ^{*} $ are such that
\begin{equation*}
\{|\tilde{k}_{i}|^{s} \theta _{i-1} \} \in \ell ^{2}, \forall s>0
\end{equation*}
\item The arguments of $ \hat{F}^{j}_{i}(k^{j}_{i}) \in \C \setminus \{ 0 \} $, $\f _{i}$, are equal up to $O(N_{i-1} ^{-\infty})$.
\end{enumerate}
Thus, the cocycles in the
K.A.M. regime are parametrized by the action of conjugations in $\VV$
composed with constant ones on the right, and the parameters
$\theta _{i} $, $a_{i} $ and $\f _{i}$, that have nonetheless to respect the
limitations of a K.A.M. scheme with given parameters.

Classification up to $H^{\s }$ conjugation is obtained by replacing
the $O(N_{i} ^{-\infty})$ by the respective $O(N_{i} ^{-\s})$ ones.
\end{thm}
Of course, the K.A.M. scheme has some tolerance with respect to the
size of some of the parameters. For example, the inequality
$N_{i-1} <  k^{j}_{i+1 }  \leq N_{i}, j=1,2$ can be violated to
a certain extent without any significant consequences, so this
classification should be taken with a grain of salt.

We observe that every representative of a constant cocycle $(\a , A_{d} )$
with $A_{d} \in DC _{\a}$ has a finite normal form, modulo the action of
conjugations in $\VV $, and are therefore
defined by a finite number of parameters in the parameter space. Therefore,
in a certain sense, the orbits of Liouvillean cocycles in the local regime, even under
$C^{\infty}$ conjugations, are bigger than the orbits of Diophantine ones,
since they have representatives in the parameter space that are not
finitely determined for any choice of parameters for the K.A.M.
scheme. Let us call $\PP$ the space of K.A.M. normal forms modulo the
action of close-to-the-identity $C^{\infty}$ conjugations. It is formed
by the data $\{ k_{i} , \epsilon _{i} , \phi _{i}  , \theta _{i}\}$,
where $k_{i}$ is the resonant frequency, $\epsilon _{i}$ is the
distance from the exact resonance, $\phi _{i} $ is the argument of
the resonant mode, and $\theta _{i}$ is the angle defined in fig.
\ref{fig KAM}. Let us also call $\tilde{\PP }$ the reduced parameter
space, where we omit the parameter $\phi _{i}$ which irrelevant for
the dynamical properties of the cocycle.
\begin{cor} \label{cor count-uncount}
The orbit of each constant cocycle $(\a , A_{d} )$ with
$A_{d} \in DC _{\a}$ has countably many representatives in $\tilde{\PP }$. The
orbit of each constant cocycle $(\a , A_{l} )$ with
$A_{l} \in \mathcal{L} _{\a}$ has uncountably
many representatives in the same space.
\end{cor}
This corollary shows in fact that the orbit of a Liouville constant exits
the K.A.M. regime and re-enters many more times than the orbit of a
Diophantine one, which is constrained by the differentiable rigidity
theorem, \ref{thm dif rig}, to leave to infinity as the norms of
conjugations grow. We also obtain the following corollary.
\begin{cor} \label{cor disc of classes}
In $\PP $, every conjugation class is dense. Every class is totally
disconnected in $\tilde{\PP}$
\end{cor}
The density could be viewed as an infinite-dimensional
analogue of the density of $x + \beta \Z \mod 1$ in $[0,1]$, for each
$x \in [0,1]$ and for $\beta \in \R \setminus \Q$ fixed, though the analogy
is quite loose.

We topologize the parameter space in a way that is compatible
with the mapping of the parameters into a space of $C^{\infty}$
functions, i.e. closeness means $O(N^{-\infty}_{n_{i}})$-closeness, and
the $h^{s}$ norms use $N_{n_{i}}^{s}$ as weights.
The formal definition would be tedious and we omit it.

The infinite dimensionality comes from the infinite number of 
significant "rotation numbers" $a_{i}$ at each step of the K.A.M. scheme.
This theorem and its corollary show that there should be no
reasonable way of defining a fibered rotation number for non-reducible
cocycles, as we can for $SL(2,\R)$ cocycles, see \cite{Herm83} and
\cite{JohMos82}.

\begin{thm} \label{thm discon of conj cl}
Every conjugacy class is dense in $\NN $. Total
disconnectedness is lost because of the action of
conjugations in $\VV $. Conjugacy classes are not locally
connected around any point.
\end{thm}
The proof of these last results implies the next theorem, which
illustrates at what point all classes and all dynamical behaviours
are indistinguishable, at least before having iterated the dynamical
system an infinite number of times.
\begin{thm} \label{thm almost conj}
Given any cocycle $(\a , A(\. )) \in \NN $ and every
conjugacy class $\mathcal{C}$ represented in $\NN $, the
cocycle $(\a , A(\. ))$ is almost conjugate to $\mathcal{C}$.
\end{thm}
The precise definition of \textit{almost conjugation}, a
generalization of almost reducibility, is given in def.
\ref{def almost conj}. The theorem is in fact slightly stronger than a corollary of the almost reducibility theorem, and we stress it since
the analysis
of the K.A.M. normal form shows that, in fact, any class of cocycles
can serve as the linear model, admittedly using as a basis the class
of constant cocycles.

All of the above theorems hold for any fixed number of frequencies
$d \in \N ^{*}$ and a bigger $d$ only results in a smaller
neighborhood of
constants where they hold true, due to Sobolev injection theorems.
If, now, we restrict ourselves to the one-frequency case ($d=1$), we can use the
powerful tool of renormalization (\cite{Krik2001}, \cite{AK2006},
\cite{FK2009}, see also \cite{NKPhD}) which we can combine with the
work of \cite{Fra2000}, \cite{Fra2004} and \cite{deAldecoa13} and
obtain the following picture, which fills the total space
$SW^{\infty}_{\a } (\T , G)$, provided that $\a \in RDC$.
\begin{thm} [\cite{Krik2001}, \cite{Fra2004}, \cite{NKPhD}]
The picture described in theorems \ref{thm spectral dichotomy} up to
\ref{thm almost conj} holds true in an open dense subset of the total space $SW^{\infty}_{\a } (\T , G)$, provided that
$\a \in RDC$ (we stress that $d=1$ in the global theorems).
\end{thm}
We have also identified the cocycles in the complementary set to
that which is renormalized into the K.A.M. regime (under the standing
arithmetic assumption).
\begin{thm} [\cite{Krik2001}, \cite{NKPhD}]
The total space $SW^{\infty}_{\a } (\T , G)$, $\a \in RDC$, is filled up by
the countable union of immersed Fr\'{e}chet manifolds corresponding
to the conjugacy classes of periodic geodesics of $G$ of
degree $r \in \N ^{*}$. These manifolds are of codimension $2r$.
\end{thm}
The spectral properties of these cocycles where studied by K.
Fraczek and R. T. de Aldecoa.
\begin{thm} [\cite{Fra2000}, \cite{Fra2004}, \cite{deAldecoa13}]
The cocycles described in the previous theorem have purely absolutely
continuous spectrum when restricted in 
$L^{2}(\T) \times \EE _{2m+1}$, for every $m \in \N$.
\end{thm}

%

Finally, we think that a picture similar to the one above should hold when
$G=SU(2)$ is replaced by any semisimple compact Lie group (cf.
\cite{KrikAst}, \cite{NKPhD}), at least in the K.A.M. regime. The
phenomena observed in the more general case should consist of
combinations and interactions between different behaviors obsverved in $SU(2)$,
respecting the conditions of linear dependence and (non-)commutativity
between the different root spaces. However, in the neighborhood of
singular geodesics interesting phenomena may appear, caused by
the interaction of the strong mixing with the weak mixing part of
the dynamics. The analysis of such systems seems to be difficult.

Further, and certainly non-exhaustive, literature in the subject includes the works
of Cl. Chavaudret (\cite{Chav11},\cite{Chav12},\cite{Chav13}), also in
collaboration with St. Marmi (\cite{ChavMarm12}) and with L. Stolovich
(\cite{ChavStol}), H. Eliasson (\cite{El1988}, \cite{El1992},
\cite{El2001}, \cite{El2002b}, see also \cite{ElAst}),
X. Hou and J. You (\cite{HY09}, \cite{HY12}) and of
G. Popov (\cite{HPop13}), Q. Zhou (\cite{YouZhou13}), and the paper of
Avila-Fayad-Kocsard \cite{AFKo2012}, which triggered this finer
study that we took up in our recent papers.

\textbf{Acknowledgment}: This work was supported by a
Capes/PNPD scholarship.
The author would like to thank Jean-Paul Thouvenot for
motivating this paper and for his limitless disposition
to explain and discuss mathematics, and Alejandro Kocsard for
the useful discussions during the preparation of the article.

\section{Notation and definitions}

\subsection{The group $SU(2)$}
The matrix group $ G = SU(2) \approx \Sp ^{3}\subset \C^{2} $ is the multiplicative group of unitary
$2 \times 2$ matrices of determinant $1$.
We will denote the matrix $S\in G$, $S=\begin{pmatrix}
z & w \\ 
-\bar{z} & \bar{w}
\end{pmatrix}
$, where $(z,w)\in \C^{2}$ and $|z|^{2}+|w|^{2}=1$, by $\{z,w\}_{G }$. The subscript will be
suppressed from the notation, unless necessary. When coordinates in $\C ^{2}$ are fixed, the
circle $\Sp ^{1}$ is naturally embedded in $G$ as the group of diagonal matrices, which is a maximal torus (i.e. a
maximal abelian subgroup) of $G$.

The Lie algebra $g=su(2)$ is naturally isomorphic to $\R^3 \approx \R \times \C$ equipped with its vector and scalar product.
The element $s=\begin{bmatrix}
it & u \\ 
-\bar{u} & -it
\end{bmatrix}
$ will be denoted by $ \{t,u\}_{g} \in \R \times \C$. The scalar product is defined by
\begin{equation*}
\langle \{t_{1} ,u_{1} \} ,\{t_{2} ,u_{2} \} \rangle =
t_{1}t_{2}+\mathcal{R} (u_{1}\bar{u}_{2})=t_{1}t_{2}+\mathcal{R} u_{1}.\mathcal{R} u_{2}+
\mathcal{I} u_{1}.\mathcal{I} u_{2}
\end{equation*}
Mappings with values in $g$ will be denoted by
\begin{equation*}
U(\. ) = \{ U_{t}(\. ), U_{z}(\. ) \}_{g}
\end{equation*}
in these coordinates, where $U_{t}(\. )$ is a real-valued and $U_{z}(\. )$ is a complex-valued function.

The adjoint action of the group on its algebra is pushed-forward to the
action of $SO(3)$ on $\R \times \C$. In particular, the diagonal matrices, of the form $S = \exp (\{2i \pi s,0\}_{g})$, we have
$Ad(S).\{ t,u\} =\{ t,e^{4i\pi s }u\} $.

\subsection{Functional Spaces}

We will consider the space $C^{\infty }(\T ,g)$ equipped with the
standard maximum norms%
\begin{equation*}
\left\Vert U\right\Vert _{s} = \max_{0\leq \sigma \leq s} \max_{\T }\left\vert \partial
^{\sigma }U(\. )\right\vert
\end{equation*}
for $s\geq 0$, and the Sobolev norms
\begin{equation*}
\left\Vert U\right\Vert _{H^{s}}^{2}=\sum_{k\in \Z^{d}}(1+|k|^{2})^{s}|\hat{U}(k)|^{2}
\end{equation*}
where $\hat{U}(k)=\int U(\. )e^{-2i\pi kx}$ are the Fourier coefficients
of $U(\. )$. The fact that the injections $H^{s +d/2}(\T ^{d},g) \hra C^{s}(\T ^{d} ,g) $
and $C^{s}(\T ^{d},g) \hra H^{s}(\T ^{d},g)$ for all $s \geq 0$ are continuous is classical.
By abusing the notation, we will note $H^{0} = L^{2}$. We will denote the
corresponding spaces of complex sequences by lowercase letters,
\begin{equation*}
h ^{s} = \{ f\in \ell ^{2}, \sum (1+n)^{2s}|f_{n}| ^{2} <
 \infty\}
\end{equation*}

For this part, see \cite{FollandHarmAn} and \cite{SteinWeissFourierEuc}. In view of the identification $G 
\approx \Sp ^{3}$, with normalized measure,
the space $C^{\infty } ( (G)  )$ of smooth $\C$-valued functions defined on $ G $, can be
identified with $C^{\infty } ( \Sp ^{3} )$, and the identification is an isometry between the $L^{2}$ spaces.

Let us give a convenient basis for $C^{\infty } ( \Sp ^{3} )$. Given a system of coordinates $(\z ,\w )$
in $\C ^{2}$, we can define an orthonormal basis for $\PP _{m}$,
the space of homogeneous polynomials of degree $m$, by $\{ \psi _{l,m} \}_{0\leq l \leq m}$ where
$ \psi _{l,m} (\z ,\w ) = \sqrt{\frac{(m+1)!}{l!(m-l)!}}  \z ^{l} \w ^{m-l}$. The group $G$ acts on $\PP _{m}$ by
\begin{equation*}
\{ z,w \}.\phi (\z ,\w ) = \phi (z \z + w \w , -\bar{w} \z + \bar{z} \w )
\end{equation*}
and the resulting representation is noted by $\p _{m}$.
For $m$ fixed, we can define the matrix coefficients relative to the basis by
$\p _{m}^{j,p} \{ z,\bar{z} , w , \bar{w}  \} \mapsto \langle \{ z,w \}.\psi _{j,m} ,\psi _{p,m} \rangle $. The matrix coefficients are
harmonic functions of $z,\bar{z},w,\bar{w}$, and are of bidegree $(m-p,p)$, i.e. they are homogeneous of degree $m-p$ in
$(z,w)$, and homogeneous of degree $p$ in $(\bar{z},\bar{w})$, and they generate the space $\EE _{m}$.
We thus obtain the decomposition
$L^{2} = \oplus _{m \in \N} \EE _{m}$

Therefore, given a system of coordinates in $\C ^{2}$, a function $f \in L^{2}( \Sp ^{3} ) $
can be written in the form
\begin{equation*}
f(z,\bar{z},w , \bar{w}) = \sum _{m \in \N } \sum_{0 \leq p \leq m} \sum _{0 \leq j \leq m } f^{m}_{j,p }
\p _{m}^{j,p} (z,\bar{z},w , \bar{w})
\end{equation*}
where $f^{m}_{j,p } \in \C $ are the Fourier coefficients.
The functions $\p _{m}^{j,p} (z,\bar{z},w , \bar{w})$ are the eigenfunctions of the Laplacian
on $\Sp ^{3}$ and consequently smooth (in fact real analytic), and they form an orthonormal basis for
$L^{2} (\Sp ^{3})$. In higher regularity, they generate a dense subspace of $C^{\infty}$.

The group $G $ acts on $C^{\infty } ( G  ) \equiv C^{\infty } (\Sp ^{3} ) $ by pullback: if $A \in G $ and
$\begin{pmatrix}
\z \\
\w
\end{pmatrix} \in \Sp ^{3} $, then, for $\phi : \Sp ^{3} \ra \C $,
\begin{equation*}
(A.\phi) \left(
\begin{pmatrix}
\z \\
\w
\end{pmatrix}
\right)
= \phi
\left( A^{*}
\begin{pmatrix}
\z \\
\w
\end{pmatrix}
\right)
\end{equation*}
If coordinates are chosen so that $A = \{ e^{2i\pi a} , 0 \} $ is diagonal, then
\begin{equation*}
A.\phi (\z,\w) = \phi( e^{-2i\pi a} \z , e^{2i\pi a} \w)
\end{equation*}
and $A$ then acts on harmonic functions by
\begin{equation*}
A.\p _{m}^{j,p} (z,\bar{z},w , \bar{w}) = e^{-2i\pi (m-2p) a } \p _{m}^{j,p} (z,\bar{z},w , \bar{w})
\end{equation*}
where $m -2p = m-p-p$ is the difference of the degree of homogeneity in $(\bar{z}, \bar{w})$
and $(z,w)$. Therefore, the harmonics in these coordinates are eigenvectors for the associated operator.
In particular, if $a $ is irrational, the eigenvectors for the eigenvalue $1$ are exactly the elements
$\p ^{j,m/2}_{m}$, $0 \leq j \leq m$.

The group of symmetries of $\C . \psi _{m/2,m} $ is exactly the
normalizer of $\TT $, the torus of matrices commuting with $A$.
We revisit the following lemma from \cite{NKInvDist}. It examines the
effect of changes of coordinates on the eigenvectors for the
eigenvalue $1$, $\psi _{m/2,m} $.
\begin{lem} \label{Legendre roots}
For a given $m >0 $ and even, $\pi _{m/2,m} (D.\psi _{m/2,m}) = 1$ iff
$D $ is in the normalizer $\NN _{\TT}$ of $\TT$. The derivative of the norm of the
projection at $D \equiv Id \in G \mod \NN _{\TT} $ is negative and
$\pi _{m/2,m} (D.\psi _{m/2,m}) < 1$ when $D \notin \NN _{\TT}$.
\end{lem}
\begin{proof}
Call $l=m/2$ and calculate the projection:
\begin{equation*}
\pi _{l,m} (\{ z,w \} _{G}.\psi _{l,m}) = \sum _{0}^{l} (-1)^{i}
\begin{pmatrix}
l \\
i
\end{pmatrix} ^{2}
|z|^{2(l-i)}|w|^{2i} \psi _{l,m} = p_{l} (|z|,|w|) \psi _{l,m}
\end{equation*}
The factor of the projection, $p_{l}$, is a Legendre polynomial in the variable $|z|^{2}$
and $|w|^{2}=1-|z|^{2}$. The conclusion follows from the properties
of Legendre polynomials.
\end{proof}

Returning to more general facts from calculus, the $C^{s }$ norms for functions in $C^{\infty } (G )$ are defined in a
classical way, and the Sobolev norms are defined by imposing a rate
of decay on $f_{m}^{j,p} $, the  coefficients of the harmonics in the expansion of $f $, $
\left\Vert f \right\Vert _{H^{s}}^{2}=\sum_{m,j,p }(1+m^{2})^{s}|f_{m}^{j,p}|^{2} $.

Finally, we will use the truncation operators for mappings $\T \ra g$:
\begin{eqnarray*}
T_{N}f(\. ) &=&\sum _{|k|\leq N}\hat{f}(k)e^{2i\pi k\. } \\
\dot{T}_{N}f(\. ) &=&T_{N}f(\. )-\hat{f}(0) \\
R_{N}f(\. ) &=&\sum  _{|k|>N}\hat{f}(k)e^{2i\pi k\. } 
\end{eqnarray*}
These operators satisfy the estimates
\begin{eqnarray} \label{truncation est}
\left\Vert T_{N}f(\. )\right\Vert _{C^{s}} &\leq
&C_{s}N \left\Vert f(\. )\right\Vert _{C^{s}} \\
\left\Vert R_{N}f(\. )\right\Vert _{C^{s}} &\leq &C_{s,s'} N^{s-s^{\prime }+2} \left\Vert
f(\. )\right\Vert _{C^{s^{\prime }}}
\end{eqnarray}


\subsection{Arithmetics, continued fraction expansion} \label{Arithmetics, continued fraction expansion}

For this section we refer the reader to \cite{KhinContFr}. Let us
introduce some notation. For $\a \in \R ^{*}$, define
$ |||\a ||| _{\Z } =dist(\a  ,\Z) = \min _{\Z } |\a - l| $,
$[\a  ]$ the integer part of $\a$, $\{\a\}$ its fractional part
and $G(\a)=\{\a^{-1}\}$, the Gauss map.

Consider $\a\in \T \setminus \Q$ fixed, and let
$\a_{n}=G^{n}(\a  )=G(\a_{n-1})$, $a_{n}=[\a_{n-1}^{-1}]$.
The following definition is classical.
\begin{defn} \label{def DC}
We will denote by $DC(\gamma ,\tau )$ the set of numbers $\a$ in $\T \setminus \Q$ such that for any $k\not=0$,
$|\a k| _{\Z }\geq \frac{\gamma ^{-1}}{|k|^{\tau }}$. Such numbers are called \textit{Diophantine}.
\end{defn}
The set $DC(\gamma ,\tau )$, for $\tau >2$ fixed and $\gamma \in \R _{+} ^{\ast}$ small is of positive
Haar measure in $\T $, and $\cup_{\gamma >0} DC(\gamma ,\tau )$ is of
full Haar measure. The numbers that do not satisfy any Diophantine
condition are called \textit{Liouvillean}. They form a residual set of $0$ Lebesgue measure.

The following definition concerns the preservation of
Diophantine properties when the algorithm of continued fractions
is applied to the number.
\begin{defn} \label{def RDC}
We will denote by $RDC(\gamma ,\tau )$ the full measure set of  \textit{recurrent Diophantine} numbers, i.e. the $\a $ in
$\T \setminus \Q$ such that $G^{n}(\a)\in DC(\gamma ,\tau )$ for infinitely many $n$.
\end{defn}
In contexts where the parameters $\gamma $ and $\tau $ are not significant, they will be omitted in the notation of both sets.

Finally, we will need to approximate the eigenvalues of matrices in
$G$ with iterates of $\a$, and thus need the following notion,
which is looser than $(\a, a) \in \T ^{d+1}$ being Diophantine.
\begin{defn} \label{def DCa}
We will denote by $DC_{\a}(\gamma ,\tau )$ the set of elements $A$ of
$G$ satisfying the following property. If
$A = D \{ e^{2i\pi a} ,0 \} D^{*}$ for some $D\in G$,
then for $k \neq 0$,
\begin{equation*}
|||a - k\a ||| \geq \frac{\gamma^{-1}}{|k|^{\tau }}
\end{equation*}
Such numbers are called \textit{Diophantine with respect to} $\a$.
\end{defn}

\subsection{Cocycles in $\T ^{d} \times SU(2)$}

\subsubsection{Definition of the dynamics}

Let $\a \in \T ^{d} \equiv \R ^{d} / \Z ^{d} $, $d \in \N ^{*}$, be an irrational rotation. If we also let $A(\. )\in C^{\infty}(\T ^{d} ,G)$, the couple
$(\a ,A(\. ))$ acts on the fibered space
$\T ^{d} \times G \ra \T ^{d} $ by
\begin{equation*}
(\a ,A(\. )).(x,S)=(x+\a ,A(x).S) , (x,S)\in \T ^{d} \times G
\end{equation*}
We will call such an action a
\textit{quasiperiodic cocycle over} $\a$ (or simply
a cocycle). The space of such actions is denoted by
$SW_{\a }^{\infty}(\T ^{d},G)$, most times abbreviated to
$SW_{\a }^{\infty}$. The number $d \in \N ^{*}$ is the number of
frequencies of the cocycle.

The space $\bigcup\nolimits_{\a \in \T ^{d} }SW_{\a }^{\infty }$ will be denoted by $SW^{\infty }$.
The space $SW_{\a }^{\infty }$ inherits the topology of $C^{\infty }(\T ^{d} ,G)$, and $SW^{\infty }$ has the standard
product topology of $\T ^{d} \times C^{\infty }(\T ^{d},G)$. We note that cocycles are
defined over more general maps and in more general contexts of regularity
and structure of the basis and fibers.

The cocycle acts on any product space $\T ^{d} \times E $, provided
that $G\curvearrowright E$, in an obvious way. The particular case
which will be important in this article is the representation of $G $
on $L^{2} (G)$, and the resulting action of the cocycle on
$L^{2}(\T ^{d} \times G)$.

The $n$-th iterate of the action is given by
\begin{eqnarray*}
(\a,A(\. ))^{n}.(x,S)&=&(n\a,A_{n}(\. )).(x,S)=(x+n\a 
,A_{n}(x).S)\\ 
&=& (x+n\a ,A(\. +(n-1)\a)...A(\. ).S)
\end{eqnarray*}
if $n>0$. Negative iterates are the inverses of positive
ones:
\begin{equation*}
(\a,A(\. ))^{-n} =((\a,A(\. ))^{n})^{-1} 
=(-n\a,A^{\ast }(\. -n\a)...A^{\ast }(\. -\a))
\end{equation*}

\subsubsection{Conjugation and reducibility}

The cocycle $(\a  ,A(\.  ))$ is called a constant cocycle if
$A(\. )=A\in G$ is a constant mapping. In that case,
the quasiperiodic product reduces to a simple product of matrices,
$ (\a  ,A)^{n}=(n\a  ,A^{n})$.

The group $C^{\infty }(\T ^{d},G)\hra SW^{\infty }(\T ^{d},G)$ acts by \textit{fibered conjugation}: Let
$B(\.  )\in C^{\infty }(\T ^{d},G)$ and $(\a  ,A(\. ))\in SW^{\infty }(\T ^{d},G)$. Then we define
\begin{eqnarray*}
Conj_{B(\.  )}.(\a  ,A(\.  )) &=&(0,B(\.  ))\circ (\a  ,A(\.  ))\circ (0,B(\.  ))^{-1} \\
&=& (\a  ,B(\.  +\a  ).A(\. ).B^{-1}(\.  ))
\end{eqnarray*}
The dynamics of $Conj_{B(\. )}.(\a  ,A(\. ))$ and $(\a ,A(\. ))$ are
essentially the same, since
\begin{equation*}
(Conj_{B(\.  )}.(\a  ,A(\.  )))^{n}=(n\a  ,B(\.  +n\a ).A_{n}(\.  ).B^{-1}(\. ))
\end{equation*}
\begin{defn}
Two cocycles $(\a ,A(\. ))$ and $(\a ,\tilde{A}(\.  ))$ in $SW^{\infty }_{\a }$ are \textit{conjugate} iff there
exists $B(\.  )\in C^{s}(\T ^{d},G)$ such that $(\a  ,\tilde{A}(\.  ))=Conj_{B(\.  )}.(\a  ,A(\.  ))$.
We will use the notation
\begin{equation*}
(\a  ,A(\.  ))\sim (\a  ,\tilde{A}(\.  ))
\end{equation*}%
to state that the two cocycles are conjugate to each other.
\end{defn}
Since constant cocycles are a class whose dynamics can be analysed, we give the following definition.
\begin{defn} \label{defn a.r.}
A cocycle will be called \textit{reducible} iff it is conjugate to a constant.
\end{defn}
In contrast with the greater part of the literature, in this article
reducible means that the transfer function is at least measurable,
whenever its regularity is not mentioned. In this article, cocycles
are always $C^{\infty}$ smooth, but the smoothness of conjugations
may vary from $H^{0} \equiv L^{2}$ to $C^{\infty}$.

Due to the fact that not all cocycles are reducible (e.g. generic cocycles in $\T \times \Sp ^{1}$
over Liouvillean rotations, but also cocycles over
Diophantine rotations, even though this result is hard to obtain, see \cite{El2002a}, \cite{Krik2001})
we also need the following concept, which has proved to be central in the study of such dynamical systems.
\begin{defn} \label{def almost reducibility}
A cocycle $(\a ,A(\. ))$ is said to be \textit{almost reducible} if there exists a sequence of conjugations
$B_{n}(\. ) \in C^{\infty}$, such that $Conj_{B_{n}(\. )}.(\a ,A(\. ))$ becomes arbitrarily close to constants
in the $C^{\infty }$ topology, i.e. iff there exists $(A_{n})$, a sequence in $G$, such that
\begin{equation*}
A_{n}^{\ast } \left(  B_{n}(\. +\a )A(\. )B_{n}^{\ast }(\. ) \right) \overset{C^{\infty } }{\ra } Id
\end{equation*}
\end{defn}
When this property is established in a K.A.M. constructive way, we can
compare the size of $F_{n} (\. ) \in C^{\infty} (\T ^{d} ,g )$, the
error term which makes this last limit into an equality, with the
rate of growth of the conjugation $B_{n}$, and obtain that
\begin{equation*}
Ad (B_{n} (\. ) ) .F_{n} (\. ) = B_{n} (\. ) .F_{n} (\. ) . B_{n}^{*} (\. ) \overset{C^{\infty } }{\ra } 0
\end{equation*}
In this case, almost reducibility in the sense of the definition above
and almost reducibility in the sense that "the cocycle can be
conjugated arbitrarily close to reducible cocycles" are equivalent.

Herein, we will prove a more general statement, concerning
conjugation close to any conjugacy class, where the same
considerations on the error term apply.
\begin{defn} \label{def almost conj}
Let $(\a ,A(\. ))$ be a given cocycle, and $\mathcal{C}$ a given
class cocycles up to conjuation. The cocycle $(\a ,A(\. ))$ is said to
be \textit{almost conjugate to $\mathcal{C}$} if there exists a sequence of conjugations
$B_{n}(\. ) \in C^{\infty}$ and a sequence of cocycles
$(\a , C_{n}(\. )) \in \mathcal{C}$, such that $Conj_{B_{n}(\. )}.(\a ,A(\. ))$ becomes arbitrarily close to $(\a , C_{n}(\. ))$ in the $C^{\infty }$
topology.
\end{defn}

\subsubsection{Review of the K.A.M. scheme and of the normal form}

\textbf{Local conjugation}. Let $(\a , Ae^{F(\. )} ) =(\a , A_{1}e^{F_{1}(\. )} ) \in SW^{\infty} (\T ,G) $ be a cocycle over a
Diophantine rotation satisfying some smallness conditions to be made more precise later on, and suppose, moreover, that
$A = \{ e^{2i\pi a} , 0 \} $ is diagonal. The goal is to conjugate the cocycle
ever closer to constant cocyles by means of an iterative scheme. This is obtained by iterating the following lemma,
for the detailed proof of which we refer to \cite{KrikAst}, \cite{El2002a} or \cite{NKPhD}. For the sake of completeness, we
sketch the proof, following the notation of \cite{NKInvDist}.
\begin{lem} \label{loc conj lem}
Let $\a  \in DC(\gamma ,\tau )$ and $K\geq C \gamma N^{\tau}$. Let, also,
$(\a  ,Ae^{F_{1}(\. )})\in SW^{\infty }(\T ^{d} ,G)$ with
\begin{equation*}
c_{0}KN ^{s_{0}} \e  _{1,0}<1
\end{equation*}
where $c_{0}, s_{0} $ depend on $\gamma , \tau $ and $d$, and
$\e  _{1,s}=\left\Vert F_{1}\right\Vert _{s}$. Then, there exists a
conjugation $G(\.  ) = G_{1}(\.  )\in C^{\infty }(\T ^{d},G)$ such that
\begin{equation*}
G_{1}(\.  +\a  ).A_{1}.e^{F_{1}(\.  )}.G_{1}^{*}(\.  )=A_{2} e^{F_{2} (\.  )}
\end{equation*}
and such that the mappings $G_{1}(\. )$ and $F_{2}(\. )$ satisfy the
following estimates
\begin{eqnarray*}
\| G_{1}(\.  )\| _{s} &\leq &  c_{1,s} (N^{s} + KN^{s+d/2}\e  _{1,0}) \\
 \e  _{2,s}&\leq &c_{2,s}K^{2}N^{2\tau +d}(N^{s}\e  _{1,0}+\e _{1,s})\e  _{1,0}+
C_{s,s^{\prime }}K^{2}N^{s-s^{\prime }+2\tau +d}\e _{1,s^{\prime }}
\end{eqnarray*}
where $s' \geq s$, and $\e _{2,s} = \| F_{2}(\.  )\| _{s}$
\end{lem}

\begin{figure}[h]
\begin{center}
\includegraphics[scale=0.33]{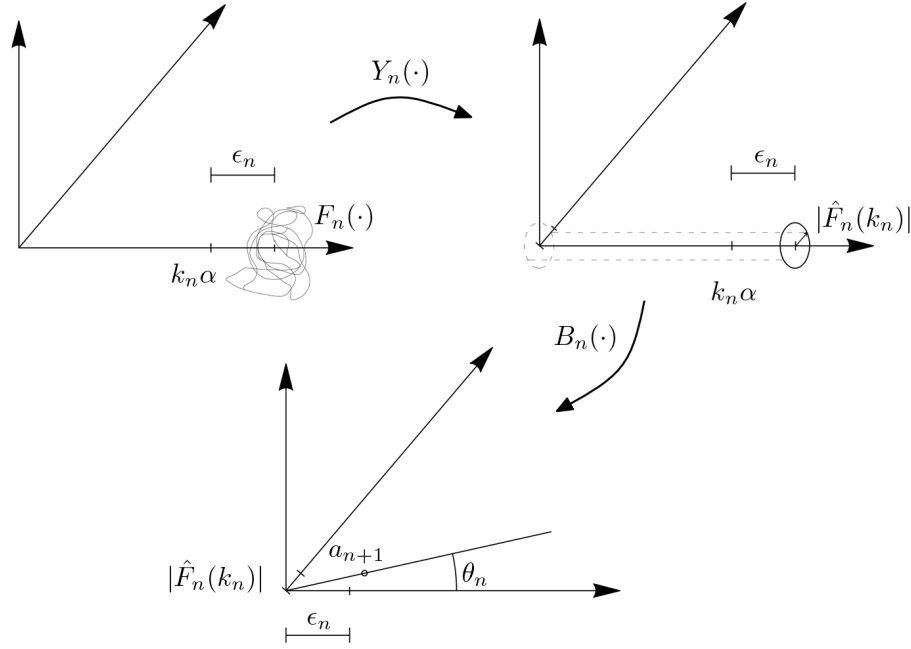} 
\end{center}
\caption{The $n$-th step of the K.A.M. scheme}
\label{fig KAM}
\end{figure}

If we suppose that $Y(\. ) : \T \ra g $ can conjugate $(\a , A_{1}e^{F_{1}(\. )} )$ to
$(\a , A_{2}e^{F_{2}(\. )} )$, with $\| F_{2}(\. ) \| \ll \| F_{1}(\. ) \| $, then it must satisfy the functional
equation
\begin{equation*}
A_{1}^{*}e^{Y(\. +\a )}A_{1}e^{F_{1}(\. )}e^{-Y(\. )}  =A_{1}^{*} A_{2}e^{F_{2}(\. )} 
\end{equation*}
Linearization of this equation under the assumption that all $C^{0}$ norms are smaller than $1$ gives
\begin{equation*} 
Ad(A_{1}^{*})Y(\. +\a ) +F_{1}(\. )-Y(\. )  = \exp ^{-1} (A_{1}^{*} A_{2})
\end{equation*}

The equation for the diagonal coordinate is a linear cohomological
one, and after truncation in the frequency domain it reads
\begin{equation*}
Y_{t } (\. +\a )-Y_{t } (\.) + \dot{T}_{N} F_{1,t }(\. )=0
\end{equation*}
and the solution as well as the estimates it satisfies are classical.
The rest satisfies the estimate of eq. \ref{truncation est}, and the
mean value $\hat{F}_{1,t } (0) $ is an obstruction and will be
integrated in $\exp ^{-1} (A_{1}^{*} A_{2}) $.

The equation concerning the non-diagonal part is a twisted
cohomological equation, whose twist depends on the linear model
$(\a , A_{1})$, and it reads
\begin{equation} \label{untruncated local eq}
e^{-4i\pi a_{1}}Y_{z }(\. +\a )-Y_{z }(\. )+F_{1,z}(\. )=0
\end{equation}
or, in the frequency domain,
\begin{equation} \label{local eq in Fourier}
(e^{2i\pi ( k\a  - 2a_{1} )}-1)\hat{Y}_{z  }(k)=-\hat{F}_{1,z }(k),~ k \in \Z
\end{equation}
If for some $k_{1}$ we have
\begin{equation*}
| k_{1 } \a  -2a_{1} | _{\Z } < K ^{-1} = N^{-\nu}
\end{equation*}
with $\nu > \tau $ to be fixed, we declare the corresponding Fourier
coefficient $\hat{F}_{1,z }(k_{1 })$ a resonance, and integrate it
to the obstructions. We know by \cite{El2002a} that
such a $k_{1 }$ (called a \textit{resonant mode}), if it exists and
satisfies $0< k_{1} \leq N $, is unique in
$\{ k \in \Z , |k- k_{1 } | \leq 2N \}$. We can thus write
\begin{equation*}
a_{1} = k_{1}\a \mod \Z + \epsilon _{1}
\end{equation*}
and call $\epsilon _{1}$
the distance to the exact resonance. If we now call $T^{k_{1 }} _{2N} $ the truncation operator
projecting on the frequencies $ 0 < |k- k_{1 } | \leq 2N$ if $k_{1 }$
exists, the equation
\begin{equation*}
e^{-4i\pi a_{1}}Y_{z }(\. +\a )-Y_{z }(\. )=-T^{k_{1 }} _{2N} F_{1,z}(\. )
\end{equation*}
can be solved and the solution satisfies the announced estimates.

In total, the equation that can be solved with good estimates is
\begin{equation*}
Ad(A_{1}^{*})Y(\. +\a ) -Y(\. ) + F_{1}(\. )= 
\{ \hat{F}_{1,t } (0), \hat{F} _{1,z} (k_{1 }) e^{2i\pi k_{1 } \.} \} +
\{ R_{N} F_{1,t} (\. ) , R^{k_{1 }} _{2N} F_{1,z} (\. ) \}
\end{equation*}
with $\| Y(\. ) \|_{s} \leq C_{s} N^{s + \nu + 1/2 } \e _{1,0} $,
and thus there exists $F_{2 }^{\prime} (\. ) $, a "quadratic" term,
such that
\begin{equation*}
e^{Y(\. +\a )}A_{1}e^{F_{1}(\. )}e^{-Y(\. )}  =
\{ e^{2i\pi (a + \hat{F}_{1,t } (0) )} ,0 \} _{G} . e^{\{ 0, \hat{F} _{1,z} (k_{1 }) e^{2i\pi k_{1 } \.} \} _{g} }
e^{F_{2}^{\prime}(\. )} 
\end{equation*}

If $k_{1 }$ exists and is non-zero, iteration of local conjugation is impossible. On the
other hand, the conjugation
$B(\. ) = \{ e^{- 2i\pi k_{1 } \. /2} ,0 \} $ is such that, if we call
$F_{1}^{\prime} (\. ) = Ad(B(\. )).F_{1}(\. ) = \{ F_{1,t}(\. ) , e^{-2i\pi k_{1} \. } F_{1,z} (\. ) \} $, similarly for $Y(\. )$,
and $A_{1}^{\prime} = B(\a ) A_{1} = \{ e^{2i\pi (a - k_{1} \a /2)} \}$,
they satisfy the equation
\begin{eqnarray*}
Ad((A')^{*})Y'(\. +\a ) -Y'(\. ) + F_{1}^{\prime}(\. ) &=& 
\{ \hat{F}_{t } (0), \hat{F} _{z} (k_{1 })  \} +
\{ R_{N} F_{t} (\. ) , e^{- 2i\pi k_{1 } \.} R^{k_{1 }} _{2N} F_{z} (\. ) \} \\
&=& \{ \hat{F}_{1,t }^{\prime} (0), \hat{F} _{1,z}^{\prime} (0)  \} +
\{ R_{N} F_{1,t}^{\prime} (\. ) , \tilde{R}^{k_{1 }} _{2N} F_{1,z} (\. ) \}
\end{eqnarray*}
where $ \tilde{R}^{k_{1 }} _{2N}$ is a dis-centered rest operator. The equation for primed variables can be
obtained from eq. \ref{untruncated local eq} by applying $Ad(B(\. ))$ and using that $B(\. )$ is a
morphism and commutes with $A_{1}$. This implies that
\begin{eqnarray*}
Conj_{B(\. )} (\a , A_{1}.\exp (\{ \hat{F}_{t } (0), \hat{F} _{z} (k_{1 }) e^{ 2i\pi k_{1 } \.}) \} ) &=&
(\a , A_{1}^{\prime}.\exp (\{ \hat{F}_{t } (0), \hat{F} _{z} (k_{1 })  \} ) \\
 &=& (\a , A_{2})
\end{eqnarray*}
that is, $B(\. )$ reduces the initial constant perturbed by the
obstructions to a cocycle close to $(\a , Id )$. If $B(\. )$ happens
to be $2$-periodic, we post-conjugate with $C(\. ) \colon 2\T \ra G$,
a torus morphism commuting with $A_{2}$ and algebraically conjugate to
\begin{equation*}
\begin{pmatrix}
e^{i\pi \.} & 0 \\
0 & e^{-i\pi \.} 
\end{pmatrix}
\end{equation*}
This conjugation adds $\pm i\pi \a $ to the arguments of the
eigenvalues of $A_{2}$ and restores $1$-periodicity without
deteriorating the estimates (see also \cite{NKPhD}), as it only
shifts the frequencies of the perturbation of $A_{2}$ by $1$. We let
the reader convince themselves that this conjugations does not
influence the estimates or the proofs of the theorems of the article
and we will omit it in the rest of the arguments.
%

\textbf{The K.A.M. scheme and normal form}. If we define the following
set of parameters, we can iterate lemma \ref{loc conj lem}. Let
$N_{n+1} = N_{n}^{1+\s } = N^{(1+\s )^{n-1}}$, where $N=N_{1}$ is big enough and
$0<\s <1$, and $K_{n} = N_{n}^{\nu }$, for some $\nu > \tau$. If we suppose that $(\a ,A_{n}e^{F_{n}(\. )})$
satisfies the hypotheses of lemma \ref{loc conj lem} for the corresponding parameters, then we obtain
a mapping $G_{n}(\. ) = B_{n}(\. ) e^{Y_{n}(\. )} $ that conjugates it to $(\a ,A_{n+1}e^{F_{n+1}(\. )})$,
and we use the notation $\e _{n,s} = \| F_{n} \|_{s}$.

If we suppose that the initial perturbation small in small norm: $\e_{1,0} < \epsilon <1$, and not big in
some bigger norm: $\e _{1, s_{0}}<1$, where $\epsilon$ and $s_{0}$ depend on the choice of parameters,
then we can prove (see \cite{NKPhD} and, through it, \cite{FK2009}), that the lemma can be iterated into a scheme,
and moreover
\begin{eqnarray*}
\e _{n,s} &=& O(N_{n}^{-\infty}) \text{ for every fixed } s \text{ and }\\
\| G_{n} \| _{s} &=& O (N_{n}^{s+\lambda}) \text{ for every } s \text{ and some fixed } \lambda >0
\end{eqnarray*}
We sum these inequalities up by saying that the norms of perturbations decay exponentially, while
conjugations grow polynomially.

This fact allows us to obtain the normal form as follows.
The product of conjugations produced by the scheme at the $n$-th step is written in the form
$H_{n}(\.  ) = B_{n}(\.  )e^{Y_{n}(\.  )}...B_{1}(\. )e^{Y_{1}(\.  )}$, where the $B_{j}(\.  )$ reduce the resonant modes.
We can rewrite the product in the form
\begin{equation*}
B_{n}(\.  )\dots B_{1}(\. ).e^{\tilde{Y}_{n}(\.  )} \cdots e^{\tilde{Y}_{2}(\.  )} e^{Y_{1}(\.  )}
\end{equation*}
where $\tilde{Y}_{j}(\.  ) = \prod _{j-1}^{1} Ad(B_{i}^{*}(\. )).  Y_{j}(\.  )$. Since the $Y_{j}(\.  )$ converge
exponentially fast to $0$ (they are conjugations comparable with $F_{j}$ with a fixed loss of derivatives)
in $C^{\infty}$, and
since the algebraic conjugation deteriorates the $C^{s}$ norms by a factor of the order of $N_{n-1}^{s+d}$,
$\prod _{\infty}^{1} \exp( \tilde{Y}_{j}(\.  ))$ always converges,
say to $D(\. ) \in C^{\infty }(\T ^{d}, G)$, even if the
$H_{n}(\. )$ do not. The cocycle $Conj_{D(\. )}(\a , A^{F(\. )})$ is the
\textit{K.A.M. normal form} of the cocycle $(\a , A^{F(\. )})$, and
it has the property that the K.A.M. scheme applied to it consists only
in the reduction of resonant modes.
\begin{notation}
For a cocycle in normal form, we relabel the indexes as
$(\a ,A_{n_{i}}e^{F_{n_{i}}}) \\ = (\a ,A_{i}e^{F_{i}})$, where $n_{i}$
is a step where a reduction of a resonant mode takes place.
\end{notation}
In the language of fig. \ref{fig KAM} a cocycle in normal form,
after the successive conjugations up to the step $i$ and in the first order
of magnitude looks like a circle around the origin in the plane tangent to
a resonant sphere $\{ S.\{2\pi( k_{i}\a + \epsilon _{i},0 \}_{g}.S^{*} \}_{S\in G} $.
Its radius is $|\hat{F}_{i}(k_{i})|$. The reduction of the resonant mode
drives $k_{i}\a$ to $0$, and reduces the rest of the perturbation
to the point of coordinates
$\{ 2\pi \epsilon _{i} , \hat{F}_{i}(k_{i}) \}_{g}$. The picture
repeats itself if we zoom in in order to see the finer scales of the
dynamics, and the first part of the picture (the reduction by
the close-to-the-identity transformation $Y_{n}(\. )$) never occurs.

At the step $i$, we will assume that the constant
$A_{i} = \{ e^{2i\pi k_{i}\a},0 \}$ is the exact resonance, and the
first order perturbation
\begin{equation*}
e^{F_{i}(\. )} = e^{\{ 2i\pi \epsilon _{i} , 0\}}. e^{\{ 0 , \hat{F}_{i}(k_{i})e^{2i\pi k_{i}\. } \}}
\end{equation*}
contains the distance from the exact resonance, $\{ 2i\pi \epsilon _{i} , 0\}$.\footnote{This choice interferes with the estimates
only when the cocycle is $C^{\infty}$ reducible.}

\section{Proof of weak mixing}

Let $f \in L^{2} (\T \times G ,\C )$ be an eigenfunction of
$U = U _{(\a , A(\. ))}$, with explicit dependence on $S$. We remark
that any eigenfunction depending non-trivially on $S$ also depends
non-trivially on $x$, unless $A(\. ) \equiv A \in G$ is constant.
Since each subspace $L^{2}(\T) \times \EE _{m}$ is invariant under $U$, for each $m$
fixed,  we can suppose that there exists $m \in \N ^{*}$ such that
$f \in  L^{2}(\T) \times \EE _{m}$.
The function $f$ then admits a development
\begin{equation*}
f(x,S) = \sum _{
\substack{
k \in \Z \\
0 \leq j,p \leq m}}
f^{m,k}_{j,p} e^{2i\pi k x} \pi ^{j,p}_{m} (z,\bar{z} , w , \bar{w}) 
= \sum _{
\substack{k \in \Z \\ 0 \leq j,p \leq m}}
f^{k}_{j,p} e^{2i\pi k x} \pi ^{j,p} (z,\bar{z} , w , \bar{w})
\end{equation*}
where we have dropped $m$ from the notation since it is considered to be fixed.
It will be replaced by $i$, the index of the step of the K.A.M. scheme.

The equation satisfied by an eigenfunction of
the Koopman operator $U$ is
\begin{equation*}
f(x-\a , A^{-1}(x).S ) = \l f(x , S)
\end{equation*}
for some fixed $\l \in \Sp ^{1}$.

For a constant cocycle $(\a , A(\. )) \equiv (\a , A)$, the following lemma
is immediate, under the assumption that
\begin{equation*}
A=
\begin{pmatrix}
e^{2i\pi a } & 0 \\
0 & e^{-2i\pi a }
\end{pmatrix}
\end{equation*}
is a diagonal matrix in the coordinates that we introduce.
\begin{lem}
Let $(\a , A)$ be a constant cocycle, and consider the canonical basis of
$\EE _{m}$, formed by the functions $\pi ^{j,p} = \pi_{m} ^{j,p}$.
Then, for every $k \in \Z$ and $0 \leq j,p \leq m$, the function
\begin{equation*}
e^{2i\pi k x } \pi ^{j,p} (z,\bar{z} , w , \bar{w}) \in  L^{2}(\T) \times \EE _{m}
\end{equation*}
is an eigenfunction of the Koopman operator $U_{(\a ,A)}$,
with eigenvalue
\begin{equation*}
e^{-2i\pi (k\a + (m-2p)a)}
\end{equation*}
\end{lem}
The proof of the lemma is by immediate calculation (or see
\cite{NKInvDist}) and points to the proof
of the main theorem of the paper, where the assumption that the
cocycle $(\a , A(\. ))$ be in K.A.M. normal form becomes relevant.
The argument is the following, and it is to be compared with the
proof of thm \ref{thm cohom stab impl red}. If $f $ is an
eigenfunction of the
operator associated to the cocycle $(\a ,A(\. ))$, and if
$(\a ,A_{i}e^{F_{i}(\. )})\overset{H^{*}_{i}} \sim (\a ,A(\. ))$ at the
$n_{i}$-th step of the K.A.M. scheme,
then, $f_{i} = f \circ (Id , H^{*}_{i})$ is an eigenfunction of
$U_{i} = U_{(\a ,A_{i}e^{F_{i}(\. )})}$.
Since $F_{i}(\. )$ is very small, $f_{i}$ should be close to an eigenfunction of the
operator $U_{i}^{\prime} = U_{(\a ,A_{i})}$. The corresponding eigenvalues
of the exact eigenfunctions will be distinct, since $\a$ is supposed Diophantine.
Since this approximation converges exponentially fast for
$n \ra \infty$, and since the support in the
frequencies in $L^{2}(\T )$ is related with the summability of the angles,
we obtain the announced theorem.

We now make the argument precise.
Clearly, $f_{0} = f$ is an eigenfunction of the operator $U_{0} = U_{(\a ,A(\. ))}$
and for the eigenvalue $\l $ iff $f_{i} = f \circ (Id , H^{*}_{i})$ is
an eigenfunction of the operator
\begin{equation*}
U_{i} = U_{(\a ,A_{i}e^{F_{i}(\. )})}, i \in \N
\end{equation*}
for the same eigenvalue. Then, linearization with respect to the dynamics gives
\begin{equation} \label{eq eigenf KAM}
\tilde{U}_{i} f _{i}  = \l f _{i} + O_{L^{2}}(N_{i}^{-\infty})
\end{equation}
where $\tilde{U}_{i} = U_{(\a ,A_{i})} $ and the constants on the $O_{L^{2}}$
depend on the norm of the function $f$. Linearization is possible
because $f$ depends in a $C^{\infty}$ (in fact real analytic) way on
the variable in $G$, and the $L^{2}$ character of the function
may only be due to the slow decay of the Fourier coefficients in
the variable in $\T $.

We will also use the fact that, if the cocycle is not measurably
reducible, resonances appear rarely.
\begin{lem}
Let $(\a ,A(\. ))$ be in normal form and not measurably reducible.
Then,
\begin{equation*}
n_{i+1} - n_{i} \ra \infty
\end{equation*}
\end{lem}
\begin{proof}
The condition that the cocycle is not measurably reducible is
equivalent to
\begin{equation*}
\{ \theta _{i} \} = \{\arctan \frac{|\hat{F}(k_{i})|}{|\epsilon _{i}|} \}
\notin \ell ^{2}
\end{equation*}
Since $|\hat{F}(k_{i})| = O(|k_{i}|^{-\infty})$ and
$|\hat{F}(k_{i})| > 0 $ for all $i $, we obtain that, also,
$|\epsilon _{i}| = O(|k_{i}|^{-\infty})$. This implies that
$||| k_{i+1} \a ||| = O(|k_{i}|^{-\infty}) $. Since, however
\begin{equation*}
\frac{\log | k_{i+1} |}{\log | k_{i} |}
\approx (1+\s )^{n_{i+1}-n_{i}}
\end{equation*}
and $\a \in DC$, we can conclude.
\end{proof}
In fact, this argument can be applied as soon as the cocycle is not
$C^{\infty}$ reducible.

Let us now apply the operator $T_{N_{i}}$ on eq.
\ref{eq eigenf KAM}
and use the fact that it commutes with $\tilde{U}_{i}$ to obtain
\begin{equation*}
\tilde{U}_{i} T_{N_{i}} f _{i}  =
\l T_{N_{i}} f _{i} + O_{L^{2}}(N_{i}^{-\infty})
\end{equation*}
Consequently, since the inverse of $\tilde{U}_{i} \circ  T_{N_{i}}$ does not
magnify the error term outside of $O_{L^{2}}(N_{i}^{-\infty})$,
$T_{N_{i}} f _{i}$ is $O_{L^{2}}(N_{i}^{-\infty})$-close
to an eigenfunction of $\tilde{U}_{i} \circ T_{N_{i}}$. Since the
eigenvalues of $\tilde{U}_{i} \circ T_{N_{i}}$ are separated by
$ \gamma ^{-1}.N^{-m\tau }_{i}$, we find that the eigenvalues of
$\tilde{U}_{i} \circ T_{N_{i}}$ and their corresponding
eigenfunctions are $O(N_{i}^{-\infty})$ and
$O_{L^{2}}(N_{i}^{-\infty})$ good approximations of the corresponding
objects for $U_{i}$. Therefore, since $L^{2}$-norms in the
original coordinates and those of the $i$-th step of the K.A.M.
scheme are the same, the same approximation holds also for the
operator $U $ and the eigenfunctions transformed accordingly.

We now compare the equations \ref{eq eigenf KAM} at the steps
$n_{i}$ and $n_{i+1} \gg n_{i} $, and examine how the divergence
of the product of conjugations sends $L^{2}$ mass to infinity,
thus contradicting the initial assumption that $f \in L^{2}$.

Let us express the eigenfunctions of $\tilde{U}_{i}$ in a
coordinate system where $(\a , A_{i})$ is diagonal. There
exists $l_{i} \in \N$, such that,  up to $O(N_{i}^{-\infty})$,\footnote{The upper index $k$ in the Fourier coefficients
$f^{i,k}_{j,p}$ is in fact redundant.}
\begin{equation*}
f_{i}(x,S_{i}) =
\sum _{
\substack{
k + (m-2p)k_{i} = l_{i} \\
0 \leq j,p \leq m \\
|k|<N_{i}
}}
f^{i,k}_{j,p} e^{2i\pi k x} \pi ^{j,p}_{i}(S_{i})
\end{equation*}
Let us, now, apply the transformation $B_{i}(\. )$ which conjugates
$(\a A_{i}e^{F_{i}(\. )})$ to $(\a A_{i+1}e^{F_{i+1}(\. )})$. In
the new coordinates,
\begin{eqnarray*}
f_{i}(x,\tilde{S}_{i}) &=&
\sum _{
\substack{
k + (m-2p)k_{i} = l_{i} \\
0 \leq j,p \leq m \\
|k|<N_{i}
}}
f^{i,k}_{j,p} e^{2i\pi (k +(m-2p)k_{i}) x} \pi ^{j,p}_{i}(\tilde{S}_{i})
\\ &=&
\sum _{
0 \leq j,p \leq m
}
\tilde{f}^{i}_{j,p} e^{2i\pi l_{i} x} \pi ^{j,p}_{i}(\tilde{S}_{i})
\end{eqnarray*}

If we let $D_{i}$ be such that $D_{i} A_{i} D_{i}^{*}$ is diagonal
in the coordinates where $A_{i+1}$ is diagonal, then $D_{i}$
is $\theta _{i}$-away from a diagonal matrix in the same
coordinates. If we also apply $D_{i}$, we obtain the new
coordinates $(x, S_{i+1})$ where the cocycle $(\a , A(\. ))$ is
represented by $(\a , A_{i+1}e^{F_{i+1}(\. )})$ and $A_{i+1}$ is
diagonal, and the expression
\begin{equation*}
f_{i}(x,S_{i+1}) =
\sum _{
\substack{
0 \leq j,p \leq m \\
|k|<N_{i}
}}
\tilde{f}^{i+1,k}_{j,p} e^{2i\pi l_{i} x} \pi ^{j,p}_{i+1}(S_{i+1})
\end{equation*}
By our observation, the formula above should coincide up to
$O_{L^{2}}(N_{i+1}^{-\infty})$ with $T_{2mN_{i}} f_{i+1}(x,S_{i+1})$,
where
\begin{equation*}
\tilde{U}_{i+1} T_{N_{i+1}} f _{i+1}  =
\l T_{N_{i+1}} f _{i+1} + O_{L^{2}}(N_{i+1}^{-\infty})
\end{equation*}
for the same eigenvalue $\l$, eventually up to $O(N_{i+1}^{-\infty})$.

The incompatibility between the two representations arises from
the transformation rule of the $\pi ^{j,p}$ under a change of basis.
More precisely, the only functions that are eigenfunctions for
the operator $\tilde{U}_{i+1}$ for the eigenvalue
$e^{2i\pi l_{i}\a}$ are the functions
$e^{2i\pi l_{i} \. } \pi ^{j,m/2}_{i+1}, 0 \leq j \leq m$, and $m$ and even number.
Therefore, the compatibility of the two expressions for the
eigenfunction would impose that
\begin{equation*}
f_{i}(x,S_{i+1}) =
\sum _{
\substack{
0 \leq j \leq m \\
|k|<N_{i}
}}
\tilde{f}^{i+1,k}_{j,m/2} e^{2i\pi l_{i} x} \pi ^{j,m/2}_{i+1}(S_{i+1})
\end{equation*}
When we insert $D_{i}^{*}$ in this expression in order to
undo the change of coordinates $\tilde{S}_{i} \mapsto S_{i+1}$,
we constrain the coefficients $f^{i,k}_{j,p}$
in the image of $\oplus _{0 \leq j \leq m} \C \pi ^{j,m/2}_{i+1}$ under
the change of coordinates, always up to $O(N_{i}^{-\infty})$.

Now, the same must hold when we compare the expressions obtained at
the steps $i+1$ and $i+2$. Comparison between the constraint on
the coefficients at the step $i+2$, i.e.
\begin{equation*}
f_{i}(x,S_{i+2}) =
\sum _{
\substack{
0 \leq j \leq m \\
|k|<N_{i}
}}
\tilde{f}^{i+2,k}_{j,m/2} e^{2i\pi l_{i} x} \pi ^{j,m/2}_{i+2}(S_{i+2})
\end{equation*}
shows that the space admissible at the step $i+1$ is restricted.
When we project the preimage of the vector
$\z _{i+2}^{m/2}\w _{i+2}^{m/2}$ under the change of coordinates
$\tilde{S}_{i+1} \mapsto S_{i+2}$, the norm of the vector shrinks
by a factor $\geq O( \theta ^{2}  _{i+1})$, i.e.
\begin{equation*}
| \langle
(D_{i+1 } )_{*} (\z _{i+2}^{m/2}\w _{i+2}^{m/2}) , \z _{i+1}^{m/2}\w _{i+1}^{m/2}
\rangle |
\leq \frac{(m/2)!}{(m+1)!}(1-O( \theta^{2} _{i+1}) )
\end{equation*}
(see lem. 2.2 of \cite{NKInvDist} and lem. \ref{Legendre roots}).

Since the different constraints on the coefficients are imposed
in different scales of the dynamics for every different $i$, or
equivalently since they correspond to
frequencies in $\Z ^{d}$ belonging to distant shells, these
constraints are independent from one another. Therefore, if the
angles are not summable in $\ell ^{2}$, the intersection of the
constraints is empty and there exists no eigenfunction in $L^{2}$.
On the other hand, if the angles are summable in $\ell ^{2}$, the
procedure converges and produces an eigenfunction as should be
expected.

\bigskip

Finally, for any given cocycle in $\NN $, we prove the
existence of a subsequence of iterates accumulating to $(0, Id)$ in the
$C^{\infty}$ topology.
\begin{prop}
All cocycles in the K.A.M. regime are rigid.
\end{prop}

\begin{proof}
Every cocycle is almost reducible to a resonant one
\begin{equation*}
(\a , A(\. )) \overset{H_{i}^{*}} \sim (\a , A_{i}) + O(N_{i+1}^{-\infty})
\end{equation*}
where $A_{i} = \{e^{2i\pi k_{i+1}\a} ,0 \} $ up to
$O(N_{i}^{-\infty})$. Since,
in the case where $(\a , A(\. ))$ is not $C^{\infty}$ reducible,
$n_{i+1} \gg n_{i}$, there exists an iterate $n_{i}< m_{i} \ll n_{i+1}$
such that $(\a , A(\. ))^{m_{i}} = ( m_{i}\a , O(N_{i}^{-\infty}) )$, and $m_{i}\a \ra 0$ when $i \ra \infty$.\footnote{In fact, $\nu > \tau $ is sufficient
for obtaining rigidity for any cocycle, but it seems cumbersome,
unless the cocycle is not $C^{\infty}$ reducible.}
\end{proof}

\section{The topology of congucacy classes}

In this section we sketch a proof of theorems
\ref{thm connect red to any in N}, \ref{thm connect red to any outside N},
\ref{thm connect any to any in N}, \ref{thm connect any to any outside N}
and \ref{thm classif of classes},
corollary \ref{cor disc of classes} and theorem \ref{thm discon of conj cl}.

The conjugations that act on a K.A.M. normal form at step $i$ of
the scheme are:
\begin{enumerate}
\item Far-from-the-identity conjugations commuting with the constant $A_{i}$
\begin{equation*}
B^{\prime}(\. )=
\begin{pmatrix}
e^{2i\pi k_{i}^{\prime}\. /2} & 0 \\
0 & e^{-2i\pi k_{i}^{\prime}\. /2}
\end{pmatrix}
\end{equation*}
where $k_{i}^{\prime}\in \Z $ is such that
that $N_{i-1} < k_{i}+k_{i}^{\prime} \leq N_{i}$.
\item Constant conjugations commuting with $A_{i}$.
\item \label{it conj activ modes} Conjugations such as those constructed in the proof of thm \ref{thm meas reduc}.
\item Conjugations regaining periodicity, if necessary.
\end{enumerate}
These conjugations act as follows.
\begin{enumerate}
\item \label{far-from-id conj on normal form} Translation of the resonance by $k_{i}^{\prime} \a /2 $ and of the
corresponding resonant frequency by $k_{i}^{\prime} $. They satisfy the
estimate $ N_{i-1}^{s+1/2} \lesssim \| B' \| _{s} \lesssim N_{i}^{s+1/2}  $.
\item Multiplication of $\hat{F}_{i}(k_{i})$ by a complex number in
$\Sp ^{1}$. This only changes the argument of
$\hat{F}_{i}(k_{i})$, and such a conjugation can be introduced as a phase in $B_{i}(\. )$. They satisfy the
estimate $\| B'' \| _{s} \lesssim N_{i-1}^{s+1/2} | \theta _{i-1}| $.
\item \label{close-to-id conj on normal form} The conjugations of the third kind are constructed as follows.
Consider a one-parameter subgroup $\{  D_{i}^{t} \}_{t\in [0,1]}$, of
minimal length such that $D_{i}A_{i+1}D_{i}^{*}$ is diagonal in the
coordinates where $A_{i}$ is diagonal and $D_{i} = D_{i}^{1}$. Then,
the path
\begin{equation*}
t \mapsto B_{i}^{*}(\. ) e^{t D_{i+1}} B_{i}(\. )
\end{equation*}
when it acts by conjugation on $(\a , A_{i}e^{F_{i}(\. )})$, transforms the
parameters of the normal form as follows
\begin{eqnarray*}
\theta _{i}^{t} &=& (1-t) \theta _{i} \\
\frac{|\hat{F} _{i}^{t}(k _{i})|}{|\epsilon _{i}^{t}|} &=& \arctan \theta _{i}^{t} \\
\sqrt{|\hat{F} _{i}^{t}(k _{i})|^{2} + (\epsilon _{i}^{t})^{2}} &\equiv &  \sqrt{|\hat{F} _{i}^{0}(k _{i})|^{2} + (\epsilon _{i}^{0})^{2}} \simeq a _{i+1}
\end{eqnarray*}
i.e. the angle between $A_{i}$ and $A_{i+1}$ is driven to zero, the
rotation number at the step $i+1$ is added to the one at the step $i$ and
the following resonances are translated by $k_{i}$.
These conjugations affect the norms by a factor
\begin{equation*}
O_{H^{s}}(| \theta ^{t} _{i} - \theta ^{0} _{i}| N_{i-1}^{s})
= O_{H^{s}}(| t \theta ^{0} _{i}| N_{i-1}^{s})
\end{equation*}
when $\theta _{i} $ is close to $0$, and therefore will be close to
the identity whenever $t$ is small enough.
\end{enumerate}
The fact that these three conjugations are the only ones who act on
the parameter space of normal forms follows from the following. In view
of item \ref{far-from-id conj on normal form}, we can assume that the
resonant mode is the same for both forms. Then, we can apply item
\ref{close-to-id conj on normal form} to each one, in order to obtain
a resonant constant, but with the resonant mode disactivated. If these
two cocycles are conjugate to each other, then their arguments
$a ^{j}_{i} + a^{j} _{i +1}$, for $j=1,2$ must be equal up to
$k^{\prime}_{i} \a$, with $k^{\prime}_{i}$ not too big (see again item
\ref{far-from-id conj on normal form}). Since resonances are unique for
this size of $k^{\prime}_{i}$, no other conjugation can act on the space
of normal forms.

These facts prove thm \ref{thm classif of classes}, by applying
the same procedure as for the construction of the K.A.M. normal form,
its corollary \ref{cor disc of classes}, and
thm \ref{thm discon of conj cl}. Corollary \ref{cor count-uncount}
is proved by combining the estimates above with the proof of thm
\ref{thm dif rig}, where it is proved that the K.A.M. scheme produces
only a finite number of resonances for cocycles reducible to a
constant in $CD_{\a }$, and, "generically", an infinite one if the
constant is in $\mathcal{L} _{\a}$.

The construction of the paths is carried out by partitioning the interval $[0,1]$ into
dyadic intervals, and then continuously deforming the parameters
of the $i$-th step of the K.A.M. scheme for $t \in [2^{i-1},2^{i}]$
in a continuous way from those of the original cocycle to those
corresponding to the normal form of the target.

Let us sketch the proof of thm \ref{thm connect red to any in N}.
First, connect the $Id$ with $A_{1}$ with a continuous path, say the
shortest one parameter connecting the two elements. This part cannot
be obtained by acting by conjugation. Then, activate corresponding
mode of the normal form by a remarametrization of, say
\begin{equation*}
t \ra \{e^{2i\pi t \epsilon _{1} }, 0 \} .
\{0 ,e^{\{0 , t\hat{F}_{1}e^{2i\pi k_{1}\. } \}} \}
\end{equation*}
Proceed by induction.

The proof of theorem \ref{thm connect red to any outside N} replaces
the first step by the following one. Let
$B_{12}^{t}(\. ) \colon [0,1/4] \ra C^{\infty} $ be a path such that
$B_{12}^{0}(\. ) \equiv Id $, and $B_{12}^{1/4}(\. ) \equiv B_{1}.B_{2}(\. )$.  This is possible, since $B_{1}.B_{2}(\. )$ is homotopic to constant
mappings. This conjugation transforms the constant cocycle $(\a ,Id)$
into $(\a , \{e^{2i\pi (a_{1}+a_{2})} ,0 \})$, where we recall that
the constants in the normal form are exactly resonant, but the path
exists $\NN$. Then, we can activate the mode in the
normal form with a path $[1/4 , 1/2] \ra C^{\infty}$ acting by a
conjugation as in item \ref{it conj activ modes}. This will drive the
constant part to $\{ e^{2i\pi k_{1}\a} , 0\} = A_{1}$, and add the
pertrurbation
\begin{equation*}
\{e^{2i\pi \epsilon _{1} } , 0 \} .
\{0 ,e^{\{0 , \hat{F}_{1}e^{2i\pi k_{1}\. } \}} \}
\end{equation*}
Proceed by induction.

The proofs of the two other similar theorems
\ref{thm connect any to any in N} and \ref{thm connect any to any outside N},
are only slightly more complicated.
%
%

\bibliography{aomsample}
\bibliographystyle{aomalpha}
\end{document}